\documentclass[fleqn]{mat01}
\usepackage{times,mathtimy,amssymb,boxit}%,subtabfg,rangecite}
\begin{document}

\setcounter{page}{333}
\firstpage{333}

\font\xxx=msam10 at 10pt
\def\ab{\mbox{\xxx{\char'245}}}

\font\zzzz=tibi at 10.4pt

\newcommand{\R}{\mathbb{R}}

\newtheorem{theo}{Theorem}
\renewcommand\thetheo{\arabic{section}.\arabic{theo}}
\newtheorem{theor}[theo]{\bf Theorem}
\newtheorem{lem}[theo]{Lemma}
\newtheorem{propo}{\rm PROPOSITION}
\newtheorem{rema}[theo]{Remark}
\newtheorem{defn}[theo]{\rm DEFINITION}
\newtheorem{exam}{Example}
\newtheorem{coro}[theo]{\rm COROLLARY}
\def\conjecture{\trivlist\item[\hskip\labelsep{\it Conjecture.}]}
\def\assump{\trivlist\item[\hskip\labelsep{\it Assumption.}]}

\renewcommand{\theequation}{\thesection\arabic{equation}}

\title{Vibrations of thin piezoelectric shallow shells: Two-dimensional
approximation}

\markboth{N~Sabu}{Vibrations of thin piezoelectric shallow shells}

\author{N~SABU}

\address{T.I.F.R. Centre, IISc Campus, Bangalore 560 012, India\\
\noindent E-mail: sabu@math.tifrbng.res.in}

\volume{113}

\mon{August}

\parts{3}

\Date{MS received 21 January 2003}

\begin{abstract}
In this paper we consider the eigenvalue problem for piezoelectric
shallow shells and we show that, as the thickness of the shell goes to
zero, the eigensolutions of the three-dimensional piezoelectric shells
converge to the eigensolutions of a two-dimensional eigenvalue problem.
\end{abstract}


\maketitle

\section{Introduction}

Lower dimensional models of shells are preferred in numerical
computations to three-dimensional models when the thickness of the
shells is `very small'. A lot of work has been done on the lower
dimensional approximation of boundary value and eigenvalue problem for
elastic plates and shells (cf.~\cite{BCM,CK,CL1,CLM,CM,SKNS1,SKNS2}). Recently some work
has been done on the lower dimensional approximation of boundary value
problem for piezoelectric shells (cf. \cite{BH1}).

In this paper, we would like to study the limiting behaviour of the
eigenvalue problems for thin piezoelectric shallow shells. We begin with
a brief description of the problem and describe the results obtained.

Let $\hat{\Omega}^\epsilon=\Phi^\epsilon(\Omega^\epsilon),
\Omega^\epsilon=\omega\times(-\epsilon, \epsilon)$ with $\omega\subset
\R^2,$ and the mapping
$\Phi^\epsilon: \overline{\Omega}^\epsilon\rightarrow\R^3$ is given by
\begin{equation*}
\Phi^\epsilon(x^\epsilon)=(x_1, x_2, \epsilon\theta(x_1, x_2))+
x_{3}^{\epsilon} a_3^\epsilon(x_1, x_2)
\end{equation*}
for all $x^\epsilon=(x_1, x_2,
x_3^\epsilon)\in\overline{\Omega}^\epsilon$, where $\theta$ is an
injective mapping of class $C^3$ and ${a}_3^\epsilon$ is a unit normal
vector to the middle surface $\Phi^\epsilon(\overline{\omega})$ of the
shell. Let $\gamma_0, \gamma_e\subset\partial\omega$ with
meas($\gamma_0)>0$ and meas($\gamma_e)>0$. Let $\hat{\Gamma}^\epsilon_0=
\Phi^{\epsilon}(\gamma_0 \times(-\epsilon, \epsilon))$ and let
$\hat{\Gamma}^\epsilon_e= \Phi^{\epsilon}(\gamma_e\times(-\epsilon,
\epsilon))$. The shell is clamped along the portion
$\hat{\Gamma}^\epsilon_0 $ of the lateral surface.

Then the variational form of the eigenvalue problem consists of finding
the displacement vector ${u}^\epsilon$, the electric potential
${\varphi}^{\epsilon}$ and $\xi^\epsilon\in\R$ satisfying
eq.~(\ref{eq:a9}). We then show that the component of the eigenvector
involving the electric potential $\varphi^\epsilon$ can be uniquely
determined in terms of the displacement vector $u^\epsilon$ and the
problem thus reduces to finding $(u^\epsilon, \xi^\epsilon)$ satisfying
equations~(\ref{eq:a31}) and (\ref{eq:a32}).

After making appropriate scalings on the data and the unknowns, we
transfer the problem to a domain $\Omega = \omega\times(-1, 1)$ which is
independent of $\epsilon$. Then we show that the scaled eigensolutions
converge to the solutions of a two-dimensional eigenvalue problem
(\ref{eq:e55}).

\section{The three-dimensional problem}

Throughout this paper, Latin indices  vary over the set $\{1,2,3\}$ and
Greek indices over the set $\{1,2\}$ for the components of vectors and
tensors. The summation over repeated indices will be used.

Let $\omega\subset \R^2$ be a bounded domain with a Lipschitz continuous
boundary $\gamma$ and let $\omega$ lie locally on one side of $\gamma$.
Let $\gamma_0, \gamma_e\subset\partial\omega$ with meas$(\gamma_0)>0$
and meas($\gamma_e)>0$. Let $\gamma_1=\partial\omega\backslash \gamma_0$
and $\gamma_s=\partial\omega\backslash \gamma_e$. For each $\epsilon >
0$, we define the sets
\begin{align*}
\Omega^{\epsilon} &= \omega\times (-\epsilon, \epsilon),\quad
\Gamma^{\pm,\epsilon}=\omega\times\{\pm\epsilon\},\quad
\Gamma^\epsilon_0=\gamma_0\times(-\epsilon, \epsilon),\\[.2pc]
\Gamma^\epsilon_1 &= \gamma_1\times(-\epsilon, \epsilon),\quad
\Gamma^\epsilon_e=\gamma_e\times(-\epsilon, \epsilon),\quad
\Gamma^\epsilon_s=\gamma_s\times(-\epsilon, \epsilon).
\end{align*}

Let $x^\epsilon=(x_1, x_2, x_3^\epsilon)$ be a generic point on
$\Omega^\epsilon$ and let $\partial_\alpha=\partial_\alpha^\epsilon=
\frac{\partial}{\partial x_\alpha}$ and
$\partial_3^\epsilon=\frac{\partial} {\partial x_3^\epsilon}$.

We assume that for each $\epsilon$, we are given a function
$\theta^\epsilon :\omega\rightarrow \R$ of class $C^3$. We then define
the map $\phi^\epsilon: \omega\rightarrow \R^3$ by
\begin{equation}
\phi^\epsilon(x_1, x_2)=(x_1, x_2, \theta^\epsilon(x_1, x_2))\quad
\mbox{for all}\ \, (x_1, x_2)\in \omega.\label{eq:a1}
\end{equation}

At each point of the surface $S^\epsilon=\phi^\epsilon(\omega)$, we
define the normal vector
\begin{equation*}
a^\epsilon= (|\partial_1\theta^\epsilon|^2+|\partial_2\theta^\epsilon|^2+1)
^{-1/2}(-\partial_1\theta^\epsilon, -\partial_2\theta^\epsilon, 1).
\end{equation*}

For each $\epsilon>0$, we define the mapping
$\Phi^\epsilon:\Omega^\epsilon \rightarrow \R^3 $ by
\begin{equation}
\Phi^\epsilon(x^\epsilon)=\phi^\epsilon(x_1, x_2)+x_3^\epsilon
a^\epsilon(x_1, x_2)\quad \mbox{for all}\ \,
x^\epsilon\in\Omega^\epsilon.\label{eq:a2}
\end{equation}

It can be shown that there exists an $\epsilon_0>0$ such that the
mappings $\Phi^\epsilon:\Omega^\epsilon\rightarrow
\Phi^\epsilon(\Omega^\epsilon)$ are $C^1$ diffeomorphisms for all
$0<\epsilon\leq \epsilon_0$. The set
$\hat{\Omega}^\epsilon=\Phi^\epsilon(\Omega^\epsilon)$ is the reference
configuration of the shell. For $0<\epsilon\leq\epsilon_0$, we define
the sets
\begin{align*}
&\hat{\Gamma}^{\pm, \epsilon} = \Phi^\epsilon(\Gamma^{\pm,\epsilon}),\quad
\hat{\Gamma}^\epsilon_0=\Phi^\epsilon(\Gamma^\epsilon_0),\quad
\hat{\Gamma}^\epsilon_1=\Phi(\Gamma^\epsilon_1),\quad
\hat{\Gamma}^\epsilon_{N} = \hat{\Gamma}^\epsilon_i\cup
\hat{\Gamma}^{\pm\epsilon},\\[.2pc]
&\hat{\Gamma}^\epsilon_e = \Phi(\Gamma^\epsilon_e),\quad
\hat{\Gamma}^\epsilon_s=\Phi(\Gamma^\epsilon_s),\quad
\hat{\Gamma}^\epsilon_{eD}=\hat{\Gamma}^\epsilon_e\cup
\hat{\Gamma}^{\pm\epsilon}
\end{align*}
and we define vectors $g^\epsilon_i$ and $g^{i,\epsilon}$ by the relations
\begin{equation*}
g^\epsilon_i=\partial^\epsilon_i\Phi^\epsilon\quad \hbox{and}\quad
g^{j,\epsilon} \cdot g^\epsilon_i=\delta^j_i
\end{equation*}
which form the covariant and contravariant basis respectively of the
tangent plane of $\Phi^\epsilon(\Omega^\epsilon)$ at
$\Phi^\epsilon(x^\epsilon)$. The covariant and contravariant metric
tensors are given respectively by
\begin{equation*}
g^\epsilon_{ij}=g^\epsilon_i \cdot g^\epsilon_j\quad \hbox{and}\quad
 g^{ij,\epsilon}=g^{i,\epsilon} \cdot g^{j,\epsilon}.
\end{equation*}

The Christoffel symbols are defined by
\begin{equation*}
\Gamma^{p,\epsilon}_{ij}=g^{p,\epsilon} \cdot \partial^\epsilon_jg^\epsilon_i.
\end{equation*}

Note however that when the set $\Omega^\epsilon$ is of the special form
$\Omega^\epsilon=\omega\times(-\epsilon, \epsilon)$ and the mapping
$\Phi^\epsilon$ is of the form (\ref{eq:a2}), the following relations
hold:
\begin{equation*}
\Gamma^{3,\epsilon}_{\alpha 3}=\Gamma^{p,\epsilon}_{33}=0.
\end{equation*}

The volume element is given by $\sqrt{g^\epsilon}{\rm d}x^\epsilon$ where
\begin{equation*}
g^\epsilon = {\rm det} (g^\epsilon_{ij}).
\end{equation*}

It can be shown that  there exist constants $g_1$ and $g_2$ such that
\begin{equation}
0 < g_1\leq g^\epsilon\leq g_2\label{eq:1a2}
\end{equation}
for $0\leq\epsilon\leq\epsilon_0$.

Let $\hat{A}^{ijkl,\epsilon}, \hat{P}^{ijk,\epsilon}$ and
$\hat{\cal E}^{ij,\epsilon}$ be the elastic, piezoelectric and dielectric
tensors respectively. We assume that the material of the shell is {\it
homogeneous and isotropic}. Then the elasticity tensor is given by
\begin{equation}
\hat{A}^{ijkl,\epsilon}=\lambda g^{ij}g^{kl} + \mu(g^{ik}g^{jl}+g^{il}g^{jk}), \label{eq:a3}
\end{equation}
where $\lambda$ and $\mu$ are the Lam\`{e} constants of the material.

These tensors satisfy the following coercive relations. There
exists a constant $C>0$ such that for all symmetric tensors $(M_{ij})$
and for any vector $(t_i)\in \R^3$,
\begin{align}
&\hat{A}^{ijkl,\epsilon}M_{kl}M_{ij} \geq C\sum_{i,j=1}^3(M_{ij})^2, \label{eq:a4}\\[.2pc]
&\hat{\cal E}^{kl,\epsilon}t_kt_l \geq C{\sum_{j=1}^3t^2_j}.\label{eq:a5}
\end{align}

Moreover we have the symmetries
\begin{equation*}
\hat{A}^{ijkl,\epsilon} = \hat{A}^{klij,\epsilon} = \hat{A}^{jikl,\epsilon},\quad
\hat{\cal E}^{kl,\epsilon} = \hat{\cal E}^{kl,\epsilon},\quad \hat{P}^{ijk,\epsilon} =
\hat{P}^{kij, \epsilon}.
\end{equation*}

Then the eigenvalue problem consists of finding
$(\hat{u}^\epsilon, \hat{\varphi}^\epsilon, \xi^\epsilon)$ such that
\begin{align}
\left. \begin{array}{lcl}
-{\rm div}\hat{\sigma}^\epsilon(\hat{u}^\epsilon, \hat{\varphi}^\epsilon)
=\xi^\epsilon\hat{u}^\epsilon \mbox{ in }\hat{\Omega}^\epsilon\\[.1pc]
\hat{\sigma}^\epsilon(\hat{u}^\epsilon, \hat{\varphi}^\epsilon)\nu
=0 \mbox{ on } \hat{\Gamma}^\epsilon_N\\[.1pc]
\hat{u}^\epsilon=0 \mbox{ on } \hat{\Gamma}^\epsilon_0
\end{array} \right\},\label{eq:aa1}\\[.2pc]
\left. \begin{array}{lcl}
{\rm div}\hat{D}^\epsilon(\hat{u}^\epsilon, \hat{\varphi}^\epsilon)=0
\mbox{ in } \hat{\Omega}^\epsilon\\[.1pc]
\hat{D}^\epsilon(\hat{u}^\epsilon, \hat{\varphi}^\epsilon)\nu=0 \mbox{
on } \hat{\Gamma}_{s}^\epsilon\\[.1pc]
\hat{\varphi}^\epsilon=0 \mbox{ on }\hat{\Gamma}^\epsilon_{eD}.
\end{array} \right\},\label{eq:aa2}
\end{align}
where
\begin{align}
\hat{\sigma}_{ij}^\epsilon &= \hat{A}^{ijkl,\epsilon}\hat{e}_{ij}^{\epsilon}
-\hat{P}^{kij,\epsilon}\hat{E}_k,\label{eq:aa3}\\[.2pc]
\hat{D}_k &= \hat{P}^{kij,\epsilon}\hat{e}_{ij}^\epsilon +
\hat{\cal E}^{kl,\epsilon}\hat{E}_l,\label{eq:aa4}
\end{align}
where $\hat{e}^\epsilon_{ij}(\hat{u}^\epsilon) =
\frac{1}{2}(\hat{\partial}^\epsilon_i\hat{u}^\epsilon_j+
\hat{\partial}^\epsilon_j\hat{u}^\epsilon_i),
\hat{\partial}^\epsilon_i = \partial/\partial\hat{x}^\epsilon_{i}$ and
$\hat{E}_k(\hat{\varphi}^\epsilon)=- \ \ \raise .01pc \hbox{$\hat{}$}\hskip -.33pc{\bigtriangledown}
(\hat{\varphi}^\epsilon)$.

We define the spaces
\begin{align}
\hat{V}^\epsilon &= \{\hat{v}\in(H^1(\hat{\Omega}^\epsilon))^3,
\hat{v}|_{\hat{\Gamma}^\epsilon_0}=0\},\label{eq:aa5}\\[.2pc]
\hat{\Psi}^\epsilon &= \{\hat{\psi}\in H^1(\hat{\Omega}^\epsilon),
\hat{\psi}|_{\hat{\Gamma}_{eD}^\epsilon}=0\}.\label{eq:aa6}
\end{align}

Then the variational form of systems (\ref{eq:aa1}) and (\ref{eq:aa2})
is to find $(\hat{u}^\epsilon, {\hat{\varphi}}^\epsilon, \xi^\epsilon)\in
\hat{V}^\epsilon\times\hat{\Psi}^\epsilon\times \R$ such that
\begin{equation}
\hat{a}^\epsilon((\hat{u}^\epsilon,
{\hat{\varphi}}^\epsilon), (\hat{v}^\epsilon, \hat{\psi}^{\epsilon}))
=\xi^\epsilon\hat{l}^\epsilon(\hat{v}^\epsilon, \hat{\psi}^\epsilon)\quad \mbox { for all }
(\hat{v}^\epsilon, \hat{\psi}^\epsilon)\in \hat{V}^\epsilon\times
\hat{\Psi}^\epsilon,\label{eq:aa7}
\end{equation}
where
\begin{align}
\hat{a}^\epsilon((\hat{u}^\epsilon, {\hat{\varphi}}^\epsilon),
(\hat{v}^\epsilon, \hat{\psi}^{\epsilon}))
&= \int_{\hat{\Omega}^\epsilon}\hat{A}^{ijkl,\epsilon}
\hat{e}_{kl}^\epsilon(\hat{u}^\epsilon)
\hat{e}_{ij}^\epsilon(\hat{v}^\epsilon){\rm d} \hat{x}^\epsilon\nonumber\\[.2pc]
&\quad\ + \int_{\hat{\Omega}^\epsilon}\hat{\cal E}^{ij,\epsilon}\hat{\partial}_i^\epsilon
{\hat{\varphi}}^\epsilon\hat{\partial}_j^\epsilon\hat{\psi}^{\epsilon}
{\rm d}\hat{x}^\epsilon\nonumber\\[.2pc]
&\quad\ +\int_{\hat{\Omega}^\epsilon} \hat{P}^{mij,\epsilon}(\hat{\partial}^\epsilon_m
{\hat{\varphi}}^\epsilon \hat{e}_{ij}^\epsilon(\hat{v}^\epsilon)-
\hat{\partial}_m^\epsilon\hat{\psi}^\epsilon
\hat{e}_{ij}^\epsilon(\hat{u}^\epsilon)){\rm d}\hat{x}^\epsilon,\label{eq:aa8}
\end{align}

$\left.\right.$\vspace{-2.5pc}

\begin{align}
\hat{l}^\epsilon(\hat{v}^\epsilon,
\hat{\psi}^{\epsilon}) &= \int_{\hat{\Omega}^\epsilon}\hat{u}^\epsilon \cdot
\hat{v}^\epsilon {\rm d}\hat{x}^\epsilon.\label{eq:aa9}
\end{align}

Since the mappings $\Phi^\epsilon:\overline{\Omega}^\epsilon\rightarrow
\overline{\hat{\Omega}}^\epsilon$ are assumed to be $C^1$ diffeomorphisms,
the correspondences that associate with every element $\hat{v}^\epsilon
\in\hat{V}^\epsilon$, the vector
\begin{equation*}
v^\epsilon=\hat{v}^\epsilon \cdot \Phi^\epsilon:
\Omega^\epsilon\rightarrow\R^3
\end{equation*}
and with every element $\hat{\psi}^{\epsilon} \in \hat{\Psi}^{\epsilon}$,
the function
\begin{equation*}
\psi^\epsilon=\hat{\psi}^\epsilon \cdot \Phi^\epsilon:\Omega^\epsilon\rightarrow\R
\end{equation*}
induce bijections between the spaces $\hat{V}^\epsilon$ and
$V^\epsilon$,  and the spaces $\hat{\Psi}^\epsilon$ and $\Psi^\epsilon$ respectively, where
\begin{align}
V^\epsilon &= \{v^\epsilon\in (H^1(\Omega^\epsilon))^3|v^\epsilon=0
\ {\rm on} \ \Gamma^\epsilon_0\},\label{eq:a6}\\[.2pc]
\Psi^\epsilon &= \{\psi^\epsilon\in H^1(\Omega^\epsilon)|\psi^{\epsilon} = 0 \ {\rm on} \
\Gamma^\epsilon_{eD}\}.\label{eq:a7}
\end{align}

Then we have
\begin{align}
&\hat{\partial}^\epsilon_j\hat{v}^\epsilon(\hat{x}^\epsilon)
= (\partial^\epsilon_i v^\epsilon)(g^{i, \epsilon})_j,\label{eq:aa11}\\[.2pc]
&\hat{e}_{ij}(\hat{v})(\hat{x}^\epsilon) = e^\epsilon_{k\|l}(v^\epsilon)
(g^{k,\epsilon})_i(g^{l,\epsilon})_j,\label{eq:aa10}
\end{align}
where
\begin{equation}
e^\epsilon_{i\|j}(v^\epsilon)=\frac{1}{2}(\partial^\epsilon_i
v^\epsilon_j+ \partial^\epsilon_j v^\epsilon_i)-\Gamma^{p,\epsilon}_{ij}v^\epsilon_p.
\label{eq:a8}
\end{equation}

Then the variational form (\ref{eq:aa7}) posed on the domain
$\Omega^\epsilon$ is to find $(u^\epsilon, {\varphi}^\epsilon, \xi^\epsilon)\in
V^\epsilon\times\Psi^\epsilon\times\R$ such that
\begin{equation}
a^\epsilon((u^\epsilon, {\varphi}^\epsilon),
(v^\epsilon, \psi^{\epsilon})) =\xi^\epsilon l^\epsilon(v^\epsilon, \psi^\epsilon )\quad \mbox { for all }
(v^\epsilon, \psi^\epsilon)\in  V^\epsilon\times\Psi^\epsilon,\label{eq:a9}
\end{equation}
where
\begin{align}
a^\epsilon((u^\epsilon, {\varphi}^\epsilon),
(v^\epsilon, \psi^{\epsilon})) &= \int_{\Omega^\epsilon}A^{ijkl,\epsilon}
e_{k\|l}^\epsilon(v^\epsilon)
e_{i\|j}^\epsilon(v^\epsilon)\sqrt{g^\epsilon}{\rm d}x^\epsilon\nonumber\\[.2pc]
&\quad\ +\int_{\Omega^\epsilon} {\cal E}^{ij,\epsilon}\partial_i^\epsilon
{\varphi}^\epsilon\partial_j^\epsilon\psi^\epsilon\sqrt{g^\epsilon}
{\rm d}x^\epsilon\nonumber\\[.2pc]
&\quad\ +\int_{\Omega^\epsilon}
P^{mij,\epsilon}(\partial^\epsilon_m{\varphi}^\epsilon
e_{i\|j}^\epsilon(v^\epsilon)\nonumber\\[.2pc]
&\quad\ -\partial_m^\epsilon\psi^\epsilon
e_{i\|j}^\epsilon(u^\epsilon))\sqrt{g^\epsilon}{\rm d}x^\epsilon,\label{eq:a10}
\end{align}

$\left.\right.$\vspace{-2pc}

\begin{align}
&l^\epsilon(v^\epsilon, \psi^{\epsilon}) = \int_{\Omega^\epsilon}u^\epsilon \cdot v^\epsilon
\sqrt{g^\epsilon} {\rm d}x^\epsilon,\label{eq:a11}\\[.2pc]
&A^{pqrs, \epsilon} = \hat{A}^{ijkl, \epsilon}(g^{p,\epsilon})_i \cdot
(g^{q,\epsilon})_j \cdot (g^{r,\epsilon})_k \cdot (g^{s,\epsilon})_l,\label{eq:aa12}\\[.2pc]
&{\cal E}^{pq,\epsilon} = \hat{\cal E}^{ij,\epsilon}(g^{p,\epsilon})_i \cdot
(g^{q,\epsilon})_j,\label{eq:aa13}\\[.2pc]
&P^{pqr,\epsilon} = \hat{P}^{ijk,\epsilon}(g^{p,\epsilon})_i \cdot
(g^{q,\epsilon})_j \cdot (g^{r,\epsilon})_k. \label{eq:aa14}
\end{align}

Using the relations (\ref{eq:1a2}), (\ref{eq:a4}) and (\ref{eq:a5}),
it can be shown that there exists a constant
$C>0$ such that for all symmetric tensor $(M_{ij})$ and for any vector
$(t_i)\in \R^3$,
\begin{align}
&A^{ijkl,\epsilon}M_{kl}M_{ij} \geq C \sum_{i,j=1}^3(M_{ij})^2,\label{eq:aa15}\\[.2pc]
&{\cal E}^{ij,\epsilon}t_it_j \geq C \sum_{i=1}^3t_i^2.\label{eq:aa16}
\end{align}

Clearly the bilinear form associated with the left-hand side of
(\ref{eq:a9}) is elliptic. Hence by Lax--Milgram theorem, given
$f^\epsilon\in V^{\prime \epsilon}$ and
$h^\epsilon\in \Psi^{\prime \epsilon}$,
there exists a unique $(u^\epsilon, \varphi^\epsilon)\in V^\epsilon\times\Psi^\epsilon$ such that
\begin{equation}
a^\epsilon((u^\epsilon, \varphi^\epsilon), (v^\epsilon,
\psi^\epsilon))= \langle (f^\epsilon, h^\epsilon), (v^\epsilon, \psi^\epsilon)\rangle\quad \ \
\forall V^{\epsilon} \times \Psi^{\epsilon} \in V^{\epsilon} \times \Psi^{\epsilon}.\label{eq:a12}
\end{equation}

In particular, for each $f^\epsilon\in (L^2(\Omega^\epsilon))^3$, there
exists a unique solution $(u^\epsilon, \varphi^\epsilon)\in
V^\epsilon\times\Psi^\epsilon$ such that
\begin{equation}
a^\epsilon((u^\epsilon, \varphi^\epsilon), (v^\epsilon,
\psi^\epsilon))=\int_{\Omega^\epsilon} f^\epsilon v^\epsilon
\sqrt{g^\epsilon}{\rm d}x^\epsilon\quad \forall v^\epsilon\times\psi^\epsilon\in
V^\epsilon\times\Psi^\epsilon. \label{eq:a13}
\end{equation}

This is equivalent to the following equations.
\begin{align}
\int_{{\Omega}^\epsilon}{A}^{ijkl,\epsilon}{e}_{k\|l}^\epsilon({u}^\epsilon)
{e}_{i\|j}^\epsilon({v}^\epsilon)\sqrt{g^\epsilon}{\rm d}x^\epsilon
&+ \int_{{\Omega}^\epsilon} {P}^{mij,\epsilon}{\partial}^\epsilon_m
(\varphi^\epsilon) {e}_{i\|j}^\epsilon({v}^\epsilon)\sqrt{g^\epsilon}{\rm d}x^\epsilon\nonumber\\
&= \int_{\Omega^\epsilon} f^\epsilon v^\epsilon
\sqrt{g^\epsilon}{\rm d}x^\epsilon\quad \forall v^\epsilon\in V^\epsilon \label{eq:a14}
\end{align}
and
\begin{equation}
\int_{{\Omega}^\epsilon}{\cal E}^{ij,\epsilon}{\partial}_i^\epsilon
{{\varphi}}^\epsilon{\partial}_j^\epsilon{\psi}^{\epsilon}
\sqrt{g^\epsilon}{\rm d}x^\epsilon=\int_{{\Omega}^\epsilon}{P}^{mij,\epsilon}
{\partial}_m^\epsilon{\psi}^\epsilon
{e}_{i\|j}^\epsilon({u}^\epsilon)\sqrt{g^\epsilon}{\rm d}x^\epsilon\ \ \ \forall
\psi^\epsilon\in \Psi^\epsilon.\label{eq:a15}
\end{equation}

From relation (2.28), it follows that
the bilinear form associated with the left-hand side of (\ref{eq:a15}) is
$\Psi^\epsilon$-elliptic.

Also for each $h^\epsilon\in V^\epsilon$, the mapping
\begin{equation*}
\psi^\epsilon\rightarrow \int_\Omega^\epsilon P^{mij,\epsilon}
\partial_m\psi^\epsilon e^\epsilon_{i\|j}(h^\epsilon)
\sqrt{g^\epsilon}{\rm d}x^\epsilon
\end{equation*}
defines a linear functional on $\Psi^\epsilon$. Hence for each
$h^\epsilon\in V^\epsilon$, there exists a unique
$T^\epsilon(h^\epsilon)\in \Psi^\epsilon$ such that
\begin{equation}
\int_{{\Omega}^\epsilon}\!\!{\cal E}^{ij,\epsilon}{\partial}_i^\epsilon
{T^\epsilon(h^\epsilon)}{\partial}_j^\epsilon{\psi}^{\epsilon}
\sqrt{g^\epsilon}{\rm d}x^\epsilon \!\!=\!\!\int_{{\Omega}^\epsilon}\!\!{P}^{mij,\epsilon}
{\partial}_m^\epsilon{\psi}^\epsilon
{e}_{i\|j}^\epsilon({h}^\epsilon)\sqrt{g^\epsilon}{\rm d}x^\epsilon\ \ \ \forall
\psi^\epsilon\!\in\! \Psi^\epsilon\label{eq:a16}
\end{equation}
and that $T^\epsilon: V^\epsilon\rightarrow \Psi^\epsilon$ is
continuous.

In particular, it follows from (\ref{eq:a15}) and the above equation that
$\varphi^\epsilon=T^\epsilon(u^\epsilon)$ and eqs~(\ref{eq:a14}) and (\ref{eq:a15}) become
\begin{align}
\int_{{\Omega}^\epsilon}\! {A}^{ijkl,\epsilon}{e}_{k\|l}^\epsilon({u}^\epsilon)
{e}_{i\|j}^\epsilon({v}^\epsilon)\sqrt{g^\epsilon}{\rm d}x^\epsilon
&+\int_{{\Omega}^\epsilon} \!{P}^{mij,\epsilon}{\partial}^\epsilon_m
(T^\epsilon(u^\epsilon)) {e}_{i\|j}^\epsilon({v}^\epsilon)\sqrt{g^\epsilon}{\rm d}x^\epsilon\nonumber\\[.2pc]
&=\int_{\Omega^\epsilon}\! f^\epsilon v^\epsilon
\sqrt{g^\epsilon}{\rm d}x^\epsilon\ \ \ \forall v^\epsilon\in V^\epsilon,\label{eq:a22}\\[.2pc]
\int_{{\Omega}^\epsilon}{\cal E}^{ij,\epsilon}{\partial}_i^\epsilon
(T^\epsilon(u^\epsilon)){\partial}_j^\epsilon{\psi}^{\epsilon}
\sqrt{g^\epsilon}{\rm d}x^\epsilon &= \int_{{\Omega}^\epsilon}{P}^{mij,\epsilon}
{\partial}_m^\epsilon{\psi}^\epsilon
{e}_{i\|j}^\epsilon({u}^\epsilon)\sqrt{g^\epsilon}{\rm d}x^\epsilon\nonumber\\[.2pc]
&\quad\ \forall \psi^\epsilon\in \Psi^\epsilon.\label{eq:a23}
\end{align}

\begin{lem}
For each $h^{\epsilon}\in (L^2(\Omega^\epsilon))^3${\rm ,}
there exists a unique $G^\epsilon(h^\epsilon) \in V^\epsilon$ such that
\begin{align}
\hskip -3.5pc \int_{{\Omega}^\epsilon}{A}^{ijkl,\epsilon}{e}_{k\|l}^\epsilon
(G^\epsilon(h^\epsilon)) {e}_{i\|j}^\epsilon({v}^\epsilon)\sqrt{g^\epsilon}{\rm d}x^\epsilon
&+ \int_{{\Omega}^\epsilon} {P}^{mij,\epsilon}{\partial}^\epsilon_m
(T^\epsilon(G^\epsilon(h^\epsilon))) {e}_{i\|j}^\epsilon({v}^\epsilon) \sqrt{g^\epsilon}{\rm d}x^\epsilon\nonumber\\[.2pc]
\hskip -3.5pc &=\int_{\Omega^\epsilon} h^\epsilon v^\epsilon
\sqrt{g^\epsilon}{\rm d}x^\epsilon\quad \forall v^\epsilon\in V^\epsilon\label{eq:a24}
\end{align}
and that $G^\epsilon:(L^2(\Omega^\epsilon))^3 \rightarrow V^\epsilon$ is continuous.
\end{lem}

\begin{proof}
Let $B^\epsilon(u^\epsilon, v^\epsilon)$ denotes the bilinear form
associated with the left-hand side of eq.~(\ref{eq:a22}). Using
(\ref{eq:a23}), we have
\begin{align}
B^\epsilon(u^\epsilon, v^\epsilon) &=
\int_{{\Omega}^\epsilon}{A}^{ijkl,\epsilon}{e}_{k\|l}^\epsilon({u}^\epsilon)
{e}_{i\|j}^\epsilon({v}^\epsilon)\sqrt{g^\epsilon}{\rm d}x^\epsilon\nonumber\\[.2pc]
&\quad\ +\int_{{\Omega}^\epsilon} {P}^{mij,\epsilon}{\partial}^\epsilon_m
(T^\epsilon(u^\epsilon)) {e}_{i\|j}^\epsilon({v}^\epsilon)\sqrt{g^\epsilon}{\rm d}x^\epsilon\nonumber\\[.2pc]
&= \int_{{\Omega}^\epsilon}{A}^{ijkl,\epsilon}{e}_{k\|l}^\epsilon({u}^\epsilon)
{e}_{i\|j}^\epsilon({v}^\epsilon)\sqrt{g^\epsilon}{\rm d}x^\epsilon\nonumber\\[.2pc]
&\quad\ +\int_{{\Omega}^\epsilon}{\cal E}^{ij,\epsilon}{\partial}_i^\epsilon
(T^\epsilon(u^\epsilon)){\partial}_j^\epsilon(T^\epsilon(v^\epsilon))
\sqrt{g^\epsilon}{\rm d}x^\epsilon\nonumber
\end{align}
\begin{align}
&= \int_{{\Omega}^\epsilon}{A}^{ijkl,\epsilon}{e}_{k\|l}^\epsilon({v}^\epsilon)
{e}_{i\|j}^\epsilon({u}^\epsilon)\sqrt{g^\epsilon}{\rm d}x^\epsilon\nonumber\\[.2pc]
&\quad\ +\int_{{\Omega}^\epsilon}{\cal E}^{ij,\epsilon}{\partial}_i^\epsilon
(T^\epsilon(v^\epsilon)){\partial}_j^\epsilon(T^\epsilon(u^\epsilon))
\sqrt{g^\epsilon}{\rm d}x^\epsilon\nonumber\\[.2pc]
&= B^\epsilon(v^\epsilon, u^\epsilon).\label{eq:a25}
\end{align}
Also, using (\ref{eq:a23}) and the relations (\ref{eq:aa15}) and (\ref{eq:aa16}),
we have
\begin{align}
B^\epsilon(u^\epsilon, u^\epsilon) &=
\int_{{\Omega}^\epsilon}{A}^{ijkl,\epsilon}{e}_{k\|l}^\epsilon({u}^\epsilon)
{e}_{i\|j}^\epsilon({u}^\epsilon)\sqrt{g^\epsilon}{\rm d}x^\epsilon\nonumber\\[.2pc]
&\quad\ +\int_{{\Omega}^\epsilon} {P}^{mij,\epsilon}{\partial}^\epsilon_m
(T^\epsilon(u^\epsilon))
{e}_{i\|j}^\epsilon({u}^\epsilon)\sqrt{g^\epsilon}{\rm d}x^\epsilon\nonumber\\[.2pc]
&= \int_{{\Omega}^\epsilon}{A}^{ijkl,\epsilon}{e}_{k\|l}^\epsilon({u}^\epsilon)
{e}_{i\|j}^\epsilon({u}^\epsilon)\sqrt{g^\epsilon}{\rm d}x^\epsilon\nonumber\\[.2pc]
&\quad\ +\int_{{\Omega}^\epsilon}{\cal E}^{ij,\epsilon}{\partial}_i^\epsilon
(T^\epsilon(u^\epsilon)){\partial}_j^\epsilon(T^\epsilon(u^\epsilon))
\sqrt{g^\epsilon}{\rm d}x^\epsilon\nonumber\\[.2pc]
&\geq C\|u^\epsilon\|^2_{V^\epsilon}.\label{eq:a26}
\end{align}
Hence $B^\epsilon(\cdots)$ is symmetric and $V^\epsilon$-elliptic.
Hence by Lax--Milgram theorem, there exists a unique
$G^\epsilon(h^\epsilon)$
satisfying (\ref{eq:a24}).
Letting $v^\epsilon=G^\epsilon(h^\epsilon)$ in (\ref{eq:a24}), we get
\begin{align}
&\int_{{\Omega}^\epsilon}{A}^{ijkl,\epsilon}{e}_{k\|l}^\epsilon
(G^\epsilon(h^\epsilon)) {e}_{i\|j}^\epsilon(G^\epsilon(h^\epsilon))\sqrt{g^\epsilon}
{\rm d}x^\epsilon\nonumber\\[.2pc]
&\quad\ + \int_{{\Omega}^\epsilon} {P}^{mij,\epsilon}{\partial}^\epsilon_m
(T^\epsilon(G^\epsilon(h^\epsilon))){e}_{i\|j}^\epsilon(G^\epsilon(h^\epsilon))
\sqrt{g^\epsilon}{\rm d}x^\epsilon\nonumber\\[.2pc]
&= \int_{\Omega^\epsilon} h^\epsilon G^\epsilon(h^\epsilon)
\sqrt{g^\epsilon}{\rm d}x^\epsilon\label{eq:a27} .
\end{align}
Using (\ref{eq:a23}), it becomes
\begin{align}
&\int_{{\Omega}^\epsilon}{A}^{ijkl,\epsilon}{e}_{k\|l}^\epsilon
(G^\epsilon(h^\epsilon)) {e}_{i\|j}^\epsilon(G^\epsilon(h^\epsilon))\sqrt{g^\epsilon}
{\rm d}x^\epsilon\nonumber\\[.2pc]
&\quad\ + \int_{{\Omega}^\epsilon}{\cal E}^{ij,\epsilon}{\partial}_i^\epsilon
(T^\epsilon(G^\epsilon(h^\epsilon))){\partial}_j^\epsilon
(T^\epsilon(G^\epsilon(h^\epsilon))) \sqrt{g^\epsilon}{\rm d}x^\epsilon\nonumber\\[.2pc]
&= \int_{\Omega^\epsilon} h^\epsilon G^\epsilon(h^\epsilon)
\sqrt{g^\epsilon}{\rm d}x^\epsilon.\label{eq:a28}
\end{align}
Using the relations (\ref{eq:aa15}) and (\ref{eq:aa16}), we have
\begin{equation}
\|G^\epsilon(h^\epsilon)\|_{V^\epsilon}^2\leq C^\epsilon
\|G^\epsilon(h^\epsilon)\|_{V^\epsilon}
\|h^\epsilon\|_{(L^2(\Omega^\epsilon))^3}.\label{eq:a29}
\end{equation}
Hence
\begin{equation}
\|G^\epsilon(h^\epsilon)\|_{V^\epsilon}
\leq C^\epsilon\|h^\epsilon\|_{(L^2(\Omega^\epsilon))^3}
\end{equation}
which implies that $G^\epsilon$ is continuous.\hfill \ab
\end{proof}

It follows from (\ref{eq:a22}) and the above lemma that
$u^\epsilon=G^\epsilon(f^\epsilon)$.
Since the inclusion $(H^1(\Omega^\epsilon))^3
\hookrightarrow (L^2(\Omega^\epsilon))^3$ is compact, it follows that
$G^\epsilon:(L^2(\Omega^\epsilon))^3\rightarrow
(L^2(\Omega^\epsilon))^3$ is compact. Also since the bilinear form
$B^\epsilon(\cdots)$ is symmetric, it follows that $G^\epsilon$ is
self-adjoint. Hence from the spectral theory of compact, self-adjoint operators, it
follows that there exists a sequence of eigenpairs
$(u^{m,\epsilon}, \xi^{m,\epsilon})_{m=1}^\infty$ such that
\begin{align}
&\int_{{\Omega}^\epsilon}{A}^{ijkl,\epsilon}
{e}_{k\|l}^\epsilon({u}^{m,\epsilon})
{e}_{i\|j}^\epsilon({v}^\epsilon)\sqrt{g^\epsilon}{\rm d}x^\epsilon\nonumber\\[.2pc]
&\quad\ + \int_{{\Omega}^\epsilon} {P}^{mij,\epsilon}{\partial}^\epsilon_m
(T^\epsilon(u^{m,\epsilon})) {e}_{i\|j}^\epsilon({v}^\epsilon)
\sqrt{g^\epsilon}{\rm d}x^\epsilon\nonumber\\[.2pc]
&= \xi^{m,\epsilon}\int_{\Omega^\epsilon} u^{m,\epsilon}
v^\epsilon \sqrt{g^\epsilon}{\rm d}x^\epsilon\quad \forall v^\epsilon\in
V^\epsilon,\label{eq:a31}\\[.2pc]
&\int_{{\Omega}^\epsilon}{\cal E}^{ij,\epsilon}{\partial}_i^\epsilon
(T^\epsilon(u^{m,\epsilon})){\partial}_j^\epsilon{\psi}^{\epsilon}
\sqrt{g^\epsilon}{\rm d}x^\epsilon\nonumber\\[.2pc]
&= \int_{{\Omega}^\epsilon}{P}^{mij,\epsilon}
{\partial}_m^\epsilon{\psi}^\epsilon
{e}_{i\|j}^\epsilon({u}^{m,\epsilon})\sqrt{g^\epsilon}{\rm d}x^\epsilon
\quad \forall \psi^\epsilon\in \Psi^\epsilon,\label{eq:a32}\\[.2pc]
&0<\xi^{1,\epsilon} \leq\xi^{2,\epsilon}\leq \cdots \leq\xi^{m,\epsilon}\leq
\cdots \rightarrow \infty,\label{eq:a33}\\[.2pc]
&\int_{\Omega^\epsilon}u^{m,\epsilon}_i
u^{n,\epsilon}_i\sqrt{g^\epsilon}{\rm d}x^\epsilon = \epsilon^3
\delta_{mn}.\label{eq:a34}
\end{align}
The sequence $\{u^{m,\epsilon}\}$ forms a complete orthonormal basis
for $(L^2(\Omega))^3$.

Define the Rayleigh quotient $R(\epsilon)(v^\epsilon)$ for
$v^\epsilon\in V^\epsilon$ by
\begin{equation}
\hskip -4pc R^\epsilon(v^\epsilon)=\frac{\int_{\Omega^\epsilon}A^{ijkl,\epsilon}
e_{k\|l}(v^\epsilon)e_{i\|j}(v^\epsilon)\sqrt{g^\epsilon}{\rm d}x^\epsilon
+\int_{{\Omega}^\epsilon} {P}^{mij,\epsilon}{\partial}^\epsilon_m
(T^\epsilon(v^\epsilon))
e_{i\|j}^\epsilon({v}^\epsilon)\sqrt{g^\epsilon}{\rm d}x^\epsilon}
{\int_{\Omega^\epsilon}v^\epsilon_i
v^\epsilon_i\sqrt{g^\epsilon}{\rm d}x^\epsilon}.\label{eq:a35}
\end{equation}

Then
\begin{equation}
\xi^{m,\epsilon}=\min_{W^\epsilon\in W^\epsilon_m}\max_{v^\epsilon\in
W^\epsilon\backslash\{0\}}R^\epsilon(v^\epsilon),\label{eq:a36}
\end{equation}
where $W_m^\epsilon$ denotes the collection of all $m$-dimensional
subspaces of $V^\epsilon$.

\section{The scaled problem}

\setcounter{equation}{0}
We now perform a change of variable so that the domain no longer
depends on $\epsilon$. With $x=(x_1, x_2, x_3)\in\Omega$, we associate
$x^\epsilon = (x_1, x_2, \epsilon x_3)\in\Omega^\epsilon$. Let
\begin{align*}
\Gamma_0 &= \gamma_0\times (-1,1),\quad \Gamma_1=\gamma_1\times
(-1,1),\quad \Gamma^{\pm}=\omega\times\{\pm 1\},\\[.2pc]
\Gamma_e &= \gamma_e\times(-1, 1),\quad \Gamma_s = \gamma_s\times(-1, 1),\\[.2pc]
\Gamma_N &= \Gamma_1\cup\Gamma^{+}\cup\Gamma^{-},\quad
\Gamma_{eD}=\Gamma^{+}\cup\Gamma^{-}\cup\Gamma_e.
\end{align*}

With the functions $\Gamma^{p,\epsilon}, g^\epsilon, A^{ijkl,\epsilon},
P^{ijk,\epsilon}, {\cal E}^{ij,\epsilon}:\Omega^\epsilon\rightarrow\R$, we
associate the functions $\Gamma^{p}(\epsilon), g^\epsilon, A^{ijkl}(\epsilon),
P^{ijk}(\epsilon), {\cal E}^{ij}(\epsilon):\Omega\rightarrow\R$ defined by
\begin{align}
\Gamma^{p}(\epsilon)(x) &:= \Gamma^{p,\epsilon}(x^\epsilon),\quad
g(\epsilon)(x)=g^\epsilon(x^\epsilon),\quad
A^{ijkl}(\epsilon)(x)=A^{ijkl,\epsilon}(x^\epsilon),\label{eq:b1}\\
P^{ijk}(\epsilon)(x) &= P^{ijk,\epsilon}(x^\epsilon),\quad
{\cal E}^{ij}(\epsilon)(x)= {\cal E}^{ij,\epsilon}(x^\epsilon).\label{eq:b2}
\end{align}

\begin{assump}
We assume that the shell is a shallow shell, i.e.
there exists a function $\theta\in C^3(\omega)$ such that
\begin{equation}
\phi^\epsilon(x_1, x_2)=(x_1, x_2, \epsilon\theta(x_1, x_2))\quad
\mbox{for all} \ \ (x_1, x_2) \in\omega,\label{eq:b6}
\end{equation}
i.e., the curvature of the shell is of the order of the thickness of the shell.

We make the following scalings on the eigensolutions.
\begin{align}
&u^{m,\epsilon}_\alpha(x^\epsilon) = \epsilon^2 u^m_\alpha(\epsilon)(x),
~~~~v_\alpha(x^\epsilon)=\epsilon^2 v_\alpha(x),\label{eq:b7}\\[.2pc]
&u^{m, \epsilon}_3(x^\epsilon) = \epsilon u^m_3(\epsilon)(x),
~~~~v_3(x^\epsilon)=\epsilon v_3(x),\label{eq:bb7}\\[.2pc]
&T^\epsilon (u^{m,\epsilon}(x^\epsilon)) =
\epsilon^3 T(\epsilon)(u^m(\epsilon)(x)),
~~~~T^\epsilon(v(x^\epsilon))=\epsilon^3 T(\epsilon)(v(x)),\label{eq:b8}\\[.2pc]
&\xi^{m,\epsilon} = \epsilon^2\xi^m(\epsilon).\label{eq:b9}
\end{align}

With the tensors $e_{i\|j}^\epsilon$, we associate the tensors
$e_{i\|j}(\epsilon)$ through the relation
\begin{equation}
e_{i\|j}^{\epsilon}(v^\epsilon)(x^\epsilon)=\epsilon^2
e_{i\|j}(\epsilon; v)(x). \label{eq:b10}
\end{equation}

We define the spaces
\begin{align}
V(\Omega) &= \{v\in (H^{1}(\Omega))^3, v|_{\Gamma_0}=0\},\label{eq:b11}\\[.2pc]
\Psi(\Omega) &= \{\psi\in H^1(\Omega), \psi|_{\Gamma_{eD}}=0\}.\label{eq:b12}
\end{align}

We denote $\varphi^m(\epsilon)=T(\epsilon)(u^m(\epsilon))$.
Then the variational equations (eqs~(\ref{eq:a31})--(\ref{eq:a34})) become
\begin{align}
&\int_\Omega A^{ijkl}(\epsilon)e_{k\|l}(\epsilon, u^m(\epsilon))
e_{i\|j}(\epsilon, v)\sqrt{g(\epsilon)}{\rm d}x\nonumber\\[.2pc]
&\quad\ +\int_\Omega P^{3kl}\partial_3{\varphi}^m
(\epsilon)e_{k\|l}(\epsilon, v)\sqrt{g(\epsilon)}{\rm d}x\nonumber\\[.2pc]
&\quad\ + \epsilon\int_\Omega P^{\alpha kl}(\epsilon)
\partial_\alpha{\varphi}^m(\epsilon)
e_{k\|l}(\epsilon, v)\sqrt{g(\epsilon)}{\rm d}x\nonumber\\[.2pc]
&= \xi^m(\epsilon)\int_\Omega [\epsilon^2 u^m_\alpha(\epsilon)v_\alpha
+u^m_3(\epsilon)v_3]\sqrt{g(\epsilon)}{\rm d}x \quad\mbox{for all }v\in
V(\Omega). \label{eq:b13}\\[.2pc]
&\int_\Omega {\cal E}^{33}(\epsilon)\partial_3{\varphi}^m(\epsilon)
\partial_3\psi \sqrt{g(\epsilon)}{\rm d}x\nonumber\\[.2pc]
&\quad\ +\epsilon\int_\Omega[{\cal E}^{3\alpha}(\epsilon)
(\partial_\alpha{\varphi}^m(\epsilon)
\partial_3\psi+\partial_3{\varphi}^m(\epsilon)\partial_\alpha\psi)]
\sqrt{g(\epsilon)}{\rm d}x\nonumber\\[.2pc]
&\quad\ +\epsilon^2\int_\Omega {\cal E}^{\alpha\beta}(\epsilon)
\partial_\alpha{\varphi}^m(\epsilon)
\partial_\beta\psi\sqrt{g(\epsilon)} {\rm d}x\nonumber\\[.2pc]
&=\int_\Omega P^{3kl}(\epsilon)\partial_3\psi e_{k\|l}(\epsilon,
u^m(\epsilon))\sqrt{g(\epsilon)}{\rm d}x \nonumber\\[.2pc]
&\quad\ +\epsilon\int_\Omega [P^{\alpha kl}(\epsilon)
\partial_\alpha\psi e_{k\|l}(\epsilon, u^m(\epsilon))]
\sqrt{g(\epsilon)}{\rm d}x \quad\mbox{for all } \psi\in \Psi(\Omega),\label{eq:b14}\\[.2pc]
&\int_\Omega [\epsilon^2u^m_\alpha(\epsilon) u^n_\alpha(\epsilon)
+u^m_3(\epsilon)u^n_3(\epsilon)]\sqrt{g(\epsilon)}{\rm d}x = \delta_{mn}.
\label{eq:b15}
\end{align}
\end{assump}

\section{Technical preliminaries}

\setcounter{equation}{0}
The following two lemmas are crucial; they play an important role in the
proof of the convergence of the scaled unknowns as $\epsilon\rightarrow
0$. In the sequel, we denote by $C_1, C_2, ..., C_n$ various constants
whose values do not depend on $\epsilon$ but may depend on $\theta$.

\setcounter{theo}{0}
\begin{lem}
The functions $e_{i\|j}(\epsilon, v) $ defined in {\rm (\ref{eq:b10})}
are of the form
\begin{align}
e_{\alpha\|\beta}(\epsilon; v) &= {\tilde{e}}_{\alpha\beta}(v)+\epsilon^2
e^{\sharp}_{\alpha\|\beta}(\epsilon; v),\label{eq:c1}\\[.2pc]
e_{\alpha\|3}(\epsilon; v) &= \frac{1}{\epsilon}\{{\tilde{e}}_{\alpha
3}(v) +\epsilon^2 e^{\sharp}_{\alpha\|3}(\epsilon; v)\},\label{eq:c2}\\[.2pc]
e_{3\|3}(\epsilon; v) &= \frac{1}{\epsilon^2}{\tilde{e}}_{33}(v),\label{eq:c3}
\end{align}
where
\begin{align}
{\tilde{e}}_{\alpha\beta}(v) &=
\frac{1}{2}(\partial_\alpha v_\beta+\partial_\beta
v_\alpha)-v_3\partial_{\alpha\beta}\theta,\label{eq:c4}\\[.2pc]
{\tilde{e}}_{\alpha 3}(v) &= {\frac{1}{2}}(\partial_\alpha v_3+\partial_3 v_\alpha),\label{eq:c5}\\[.2pc]
{\tilde{e}}_{33}(v) &= \partial_3 v_3\label{eq:c6}
\end{align}
and there exists constant $C_1$ such that
\begin{equation}
\sup_{0<\epsilon\leq\epsilon_0}\max_{\alpha,j}
\|e^{\sharp}_{\alpha,j}(\epsilon; v)
\|_{0,\Omega}\leq C_1\|v\|_{1,\Omega} \quad {\rm for\ all}\ \ v\in
V.\label{eq:c7}
\end{equation}
Also there exist constants $C_2, C_3$ and $C_4$ such that
\begin{align}
&\sup_{0<\epsilon\leq \epsilon_0}\max_{x\in\Omega}|g(x)-1|\leq
C_2\epsilon^2,\label{eq:c8}\\[.2pc]
&\sup_{0<\epsilon\leq \epsilon_0}\max_{x\in\Omega}|A^{ijkl}(\epsilon)-A^{ijkl}|
\leq C_3\epsilon^2,\label{eq:c9}
\end{align}
where
\begin{equation}
A^{ijkl}=\lambda\delta^{ij}\delta^{kl}+\mu(\delta^{ik}\delta^{jl}+
\delta^{il}\delta^{jk})\label{eq:c10}
\end{equation}
and
\begin{equation}
A^{ijkl}M_{kl}M_{ij}\geq C_4M_{ij}M_{ij}\label{eq:c11}
\end{equation}
for $0<\epsilon\leq\epsilon_0$ and for all symmetric tensors
$(M_{ij}).$
\end{lem}

\begin{proof}
The proof is based on Lemma~4.1 of \cite{BCM}.\hfill \ab
\end{proof}

From relation (\ref{eq:a5}) and definition (\ref{eq:b2}),
it follows that there exists a constant $C_5$ such that for any vector $(t_i)\in \R^3$,
\begin{equation}
{\cal E}^{ij}(\epsilon)t_it_j\geq C_5\sum_{j=1}^3t_j^2.\label{eq:b3}
\end{equation}

We assume that there exists functions $P^{kij}$ and ${\cal E}^{ij}$ such that
\begin{align}
&\sup_{0<\epsilon\leq
\epsilon_0}\max_{x\in\Omega}|P^{kij}(\epsilon)-P^{kij}|
\leq C_6\epsilon,\label{eq:b4}\\[.6pc]
&\sup_{0<\epsilon\leq \epsilon_0}\max_{x\in\Omega}|{\cal E}^{ij}
(\epsilon)- {\cal E}^{ij}| \leq C_7\epsilon.\label{eq:b5}
\end{align}

\begin{lem}
Let $\theta\in C^3(\omega)$ be a given function and let the functions
$\tilde{e}_{ij}$ be defined as in {\rm (\ref{eq:c4})--(\ref{eq:c6})}. Then
there exists a constant $C_8$ such that the following generalised Korn's
inequality holds{\rm :}
\begin{equation}
\|v\|_{1,\Omega}\leq C_8\left\lbrace \sum_{i,j}\|\tilde{e}_{ij}(v)\|^2_{0,
\Omega}\right\rbrace^{1/2}\label{eq:c12}
\end{equation}
for all $v\in V(\Omega)$ where $V(\Omega)$ is the space defined in
{\rm (\ref{eq:b11})}.
\end{lem}

\begin{proof}
The proof is based on Lemma~4.2 of \cite{BCM}.\hfill \ab
\end{proof}

\section{{\zzzz A priori} estimates}

\setcounter{equation}{0}

In this section, we show that for each positive integer $m$, the scaled
eigenvalues $\{\xi^m(\epsilon)\}$ are bounded uniformly with respect to
$\epsilon$.

Let $\varphi\in H^2_0(\omega)$. Then
\begin{equation}
v_\varphi:=(-x_3\partial_1\varphi, -x_3\partial_2\varphi, \varphi)\in
V(\Omega)\label{eq:d1}
\end{equation}

and
\begin{equation}
\tilde{e}_{\alpha\beta}(v_\varphi)=-x_3\partial_{\alpha\beta}\varphi-\varphi
\partial_{\alpha\beta}\theta, ~~~~
\tilde{e}_{i3}(v_\varphi)=0.\label{eq:d2}
\end{equation}

Hence
\begin{align}
e_{\alpha\|\beta}(\epsilon,
v_\varphi) &= -x_3\partial_{\alpha\beta}\varphi
-\varphi\partial_{\alpha\beta}\theta+O(\epsilon^2),\label{eq:d3}\\[.4pc]
e_{\alpha\|3}(\epsilon, v_\varphi) &= O(\epsilon),\label{eq:d4}\\[.4pc]
e_{3\|3}(\epsilon, v_\varphi) &= 0.\label{eq:d5}
\end{align}

We need the following lemma to prove the boundedness of the scaled eigenvalues.

\setcounter{theo}{0}
\begin{lem} There exists a constant $C_9>0$ such that
\begin{align}
&|\partial_3(T(\epsilon)(v_\varphi))|_{0,\Omega}\leq C_9
|\varphi|_{2,\omega},\label{eq:d6}\\[.4pc]
&|\epsilon\partial_{\alpha}(T(\epsilon)(v_\varphi))|_{0,\Omega}
\leq C_9 |\varphi|_{2,\omega}.\label{eq:d7}
\end{align}
\end{lem}

\begin{proof}
With the scalings (\ref{eq:b6})--(\ref{eq:b9}), the
variational equation (eq.~(\ref{eq:a16})) posed on the domain $\Omega$
reads as follows:

For each $h\in (H^1(\Omega))^3$, there exists a unique solution
$T(\epsilon)(h)\in(H^1(\Omega))^3$ such that
\begin{align}
&\int_\Omega {\cal E}^{33}(\epsilon)\partial_3T(\epsilon)(h)\partial_3\psi
\sqrt{g(\epsilon)}{\rm d}x\nonumber\\
&\quad\ +{\epsilon}\int_\Omega [{\cal E}^{\alpha 3}(\epsilon)
(\partial_\alpha T(\epsilon)(h)
\partial_3\psi+\partial_3T(\epsilon)(h)\partial_\alpha \psi)]
\sqrt{g(\epsilon)}{\rm d}x\nonumber\\
&\quad\ +\epsilon^2\int_\Omega {\cal E}^{\alpha\beta}(\epsilon)\partial_\alpha
T(\epsilon)(h) \partial_\beta\psi \sqrt{g(\epsilon)}{\rm d}x\nonumber\\
&=\int_\Omega P^{3kl}(\epsilon)
\partial_3\psi e_{k\|l}(\epsilon, h)\sqrt{g(\epsilon)}{\rm d}x\nonumber\\
&\quad\ +\epsilon\int_\Omega P^{\alpha kl}(\epsilon)
\partial_\alpha\psi e_{k\|l}(\epsilon, h)\sqrt{g(\epsilon)}{\rm d}x\quad
\forall \psi\in\Psi.
\end{align}

Taking $h=v_\varphi$ and $\psi=T(\epsilon)(v_\varphi)$ in
the above equation, we have
\begin{align}
&\int_\Omega {\cal E}^{33}(\epsilon)\partial_3T(\epsilon)(v_\varphi)
\partial_3 T(\epsilon)(v_\varphi)
\sqrt{g(\epsilon)}{\rm d}x\nonumber\\
&\quad\ +{\epsilon}\int_\Omega [{\cal E}^{\alpha 3}(\epsilon)
(\partial_\alpha T(\epsilon)(v_\varphi)
\partial_3 T(\epsilon)(v_\varphi)\nonumber\\
&\quad\ +\partial_3T(\epsilon)(v_\varphi)\partial_\alpha
T(\epsilon)(v_\varphi))] \sqrt{g(\epsilon)}{\rm d}x\nonumber\\
&\quad\ +\epsilon^2\int_\Omega {\cal E}^{\alpha\beta}(\epsilon)\partial_\alpha
T(\epsilon)(v_\varphi)
\partial_\beta T(\epsilon)(v_\varphi) \sqrt{g(\epsilon)}{\rm d}x\nonumber\\
&=\int_\Omega P^{3kl}(\epsilon)
\partial_3 T(\epsilon)(v_\varphi)
e_{k\|l}(\epsilon, v_\varphi)\sqrt{g(\epsilon)}{\rm d}x\nonumber\\
&\quad\ +\epsilon\int_\Omega P^{\alpha kl}(\epsilon)
\partial_\alpha T(\epsilon)(v_\varphi)e_{k\|l}(\epsilon, v_\varphi)
\sqrt{g(\epsilon)}{\rm d}x.
\end{align}

Using the relations (\ref{eq:b3}) and (\ref{eq:d2})--(\ref{eq:d5}),
it follows that there exists a constant $C_9>0$ such that
\begin{align}
&|\partial_{3}(T(\epsilon)(v_\varphi))|^2_{0,\Omega}+
|\epsilon\partial_{\alpha}(T(\epsilon)(v_\varphi))|^2_{0,\Omega}\nonumber\\
&\hskip 1cm \leq C_9\{|\partial_3T(\epsilon)(v_\varphi)|_{0,\Omega}
|\varphi|_{2,\omega}+
|\epsilon\partial_\alpha T(\epsilon)(v_\varphi)|_{0,\Omega}
|\varphi|_{2,\omega}\}\label{eq:d10}
\end{align}
and hence the result follows.\hfill \ab
\end{proof}

\begin{theor}[\!]
For each positive integer $m${\rm ,} there exists a constant $C(m)>0$ such
that
\begin{equation}
\xi^m(\epsilon)\leq C(m).\label{eq:d11}
\end{equation}
\end{theor}

\begin{proof}
Since problem (\ref{eq:b13}) was derived from (\ref{eq:a31}) after a
change of scale, we still have the variational characterization of the
scaled eigenvalues $\xi^m(\epsilon)$. Let $V_m$ denote the collection of
all $m$-dimensional subspaces of $V(\Omega)$. Then
\begin{equation}
\xi^m(\epsilon)=\min_{W\in V_m}\max_{v\in W}\frac{N(\epsilon)(v, v)}
{D(\epsilon)(v, v)},\label{eq:d12}
\end{equation}
where
\begin{align}
N(\epsilon)(v, v) &= \int_\Omega A^{ijkl}e_{k\|l}(\epsilon, v)
e_{i\|j}(\epsilon, v)\sqrt{g(\epsilon)}{\rm d}x\nonumber\\[.2pc]
&\quad\ +\int_\Omega P^{3kl}\partial_3
T(\epsilon)(v) e_{k\|l}(\epsilon, v)\sqrt{g(\epsilon)}{\rm d}x\nonumber\\[.2pc]
&\quad + \epsilon\int_\Omega P^{\alpha kl}
\partial_\alpha T(\epsilon)(v)e_{k\|l}(\epsilon,
v)\sqrt{g(\epsilon)}{\rm d}x,\label{eq:d13}\\[.2pc]
D(\epsilon)(v, v) &= \int_\Omega[\epsilon^2 v_\alpha v_\alpha+
v_3 v_3]\sqrt{g(\epsilon)}{\rm d}x.\label{eq:d14}
\end{align}

Let $W_m$ be the collection of all $m$-dimensional subspaces of
$H^2_0(\omega)$. Let $W\in W_m$. Define
\begin{equation}
{\bf W} = \{v_\varphi|\varphi\in W\}.\label{eq:d15}
\end{equation}
It follows that ${\bf{W}}\in V_m$. Hence, it follows from
(\ref{eq:d12}) that
\begin{equation}
\xi^m(\epsilon)\leq\min_{W\in W_m}\max_{\varphi\in W}
\frac{N(\epsilon)(v_\varphi, v_\varphi)}
{D(\epsilon)(v_\varphi, v_\varphi)}.\label{eq:d16}
\end{equation}
Now,
\begin{align}
D(\epsilon)(v_\varphi, v_\varphi) &= \int_\Omega [\epsilon^2 x_3^2
|\partial_\alpha\varphi|^2+|\varphi|^2]\sqrt{g(\epsilon)}{\rm d}x.\nonumber\\[.2pc]
&\geq \int_\omega \varphi^2{\rm d}\omega.\label{eq:d17}
\end{align}
Using the relations (\ref{eq:d3})--(\ref{eq:d5}) and Lemma~5.1, it
follows that
\begin{align}
&\int_\Omega A^{ijkl}e_{k\|l}(\epsilon, v_\varphi)
e_{i\|j}(\epsilon, v_\varphi)\sqrt{g(\epsilon)}{\rm d}x\leq
C\int_\omega |\triangle \varphi|^2 {\rm d}\omega,\label{eq:d18}\\[.2pc]
&\int_\Omega P^{3kl}\partial_3 T(\epsilon)(v_\varphi)e_{k\|l}(\epsilon,
v_\varphi)\sqrt{g(\epsilon)}{\rm d}x\leq
C\int_\omega |\triangle \varphi|^2 {\rm d}\omega,\label{eq:dd1}\\[.2pc]
&\epsilon\int_\Omega P^{\alpha kl}
\partial_\alpha T(\epsilon)(v_\varphi)e_{k\|l}(\epsilon, v_\varphi)
\sqrt{g(\epsilon)}{\rm d}x\leq
C\int_\omega |\triangle \varphi|^2 {\rm d}\omega.\label{eq:dd2}
\end{align}

Hence
\begin{align}
\xi^{m}(\epsilon) &\leq C\min_{W\in W_m}\max_{\varphi\in W}
\frac{\int_\omega |\triangle \varphi|^2 {\rm d}\omega}
{\int_\omega \varphi^2 {\rm d}\omega}\nonumber\\[.2pc]
&\leq C \lambda^m,\label{eq:d19}
\end{align}
where $\lambda^m$ is the $m$th eigenvalue of the two-dimensional
elliptic eigenvalue problem
\begin{align}
&\triangle^2 u = \lambda u\quad\mbox{ in } \omega \nonumber\\
&u=\partial_\nu u = 0 \quad\mbox{ on }\partial\omega .\label{eq:d20}
\end{align}
This completes the proof of the theorem on setting $C(m)=C\lambda^m$.\hfill \ab
\end{proof}

\section{The limit problem}

\setcounter{equation}{0}
\setcounter{theo}{0}

\begin{theor}[\!] {\rm (a)} For each positive integer $m${\rm ,} there
exists $u^m\in H^1(\Omega), \varphi^m\in L^2(\Omega)$ and $\xi^m\in\R$
such that
\begin{align}
&u^m(\epsilon)\rightarrow u^m \mbox{ in } H^1(\Omega),~~~~ \varphi^m(\epsilon)
\rightarrow \varphi^m \ {\rm in}\ L^2(\Omega),\label{eq:e1}\\[.2pc]
&(\epsilon\partial_1\varphi^m(\epsilon),
\epsilon\partial_2\varphi^m(\epsilon),  \partial_3\varphi^m(\epsilon))
\rightarrow (0, 0, \partial_3\varphi^m)\ {\rm in}\
L^2(\Omega),\label{eq:e2}\\[.2pc]
&\xi^m(\epsilon)\rightarrow \xi^m.\label{eq:e3}
\end{align}

{\rm (b)} Define the spaces
\begin{align}
&V_H(\omega) = \{(\eta_\alpha)\in (H^1(\omega))^2; \eta_\alpha=0 \ {\rm on}\
\gamma_0\},\label{eq:e4}\\[.4pc]
&V_3(\omega) = \{\eta_3\in H^2(\omega);
\eta_3=\partial_\nu\eta_3=0 \ {\rm on}\  \gamma_0 \},\label{eq:e5}\\[.4pc]
&V_{KL} = \{v\in H^1(\Omega)| v=\eta_\alpha-x_3\partial_\alpha\eta_3,
(\eta_i)\in V_H(\omega)\times V_3(\omega)\},\label{eq:e6}\\[.4pc]
&\Psi_l = \{\psi\in L^2(\Omega), \partial_3\psi\in
L^2(\Omega)\},\label{eq:e7}\\[.4pc]
&\Psi_{l0} = \{\psi\in L^2(\Omega), \partial_3\psi\in L^2(\Omega),
\psi|\Gamma^{\pm}=0\}.\label{eq:e8}
\end{align}

Then there exists $(\zeta^m_\alpha, \zeta^m_3)\in
V_H\times V_3(\omega)$ such that
\begin{align}
u^m_\alpha &= \zeta^m_\alpha-x_3\partial_\alpha\zeta^m_3 \quad {\rm and}\quad
u^m_3=\zeta^m_3,\label{eq:e9}\\[.4pc]
\varphi^m &= (1-x_3^2)\frac{p^{3\alpha\beta}}{p^{33}}
\partial_{\alpha\beta}\xi^m_3\label{eq:e10}
\end{align}
and $(\zeta^m, \xi^m)\in V_H\times V_3\times\R$ satisfies
\begin{align}
 &\hskip -4pc -\int_\omega m_{\alpha\beta}(\zeta^m)\partial_{\alpha\beta}\eta_3
{\rm d}\omega + \int_\omega
n_{\alpha\beta}^\theta(\zeta^m)\partial_{\alpha\beta}\theta\eta_3
{\rm d}\omega%\nonumber\\[.2pc]&
 + \frac{2}{3}\int_\omega \frac{p^{3\alpha\beta}p^{3\rho\tau}}{p^{33}}
\partial_{\rho\tau}\zeta^m_3\partial_{\alpha\beta}\eta_3
{\rm d}\omega \nonumber\\[.4pc]
 & = \xi^m\int_\omega \zeta^m_3\eta_3 {\rm d}\omega\quad \forall\eta_3\in V_3(\omega),\hskip 4pc
\label{eq:e11}\\[.4pc]
 & \int_\omega n_{\alpha\beta}^\theta\partial_\beta\eta_\alpha {\rm d}\omega
= 0\quad \forall\eta_\alpha\in V_H(\omega),\label{eq:e12}\hskip 4pc
\end{align}
where
\begin{align}
&m_{\alpha\beta}(\zeta) =
-\left\lbrace \frac{4\lambda\mu}{3(\lambda+4\mu)}\triangle\zeta_3\delta_{\alpha\beta}
+\frac{4\mu}{3}\partial_{\alpha\beta}\zeta_3\right\rbrace\label{eq:e13}\\[.2pc]
&n_{\alpha\beta}^\theta(\zeta) =
\frac{4\lambda\mu}{\lambda+2\mu}\tilde{e}_{\sigma\sigma}(\zeta)
\delta_{\alpha\beta}+4\mu\tilde{e}_{\alpha\beta}(\zeta)\label{eq:e14}\\[.2pc]
&p^{33} = \frac{1}{\mu}P^{3\alpha 3}P^{3\alpha 3}+\frac{1}{\lambda+2\mu}
P^{333}P^{333}+ {\cal E}^{33}\label{eq:e15}\\[.2pc]
&p^{3\alpha\beta} = P^{3\alpha\beta}-\frac{\lambda}{\lambda+2\mu}P^{333}
\delta^{\alpha\beta}.\label{eq:e16}
\end{align}
\end{theor}

\begin{proof} For the sake of clarity, the proof is divided into
several steps.\vspace{.5pc}

\noindent {\it Step} (i).\ \ Define the vector ${\tilde{\varphi}^m}_{i}(\epsilon)$ and the tensor
$\tilde{K}^m(\epsilon)=(\tilde{K}^m_{ij}(\epsilon))$ by
\begin{align}
&{\tilde{\varphi}^m}_{i}(\epsilon) =
(\epsilon\partial_1{\varphi}^m(\epsilon),
\epsilon\partial_2{\varphi}^m(\epsilon),
\partial_3{\varphi}^m(\epsilon)),\label{eq:e19}\\[.2pc]
&\tilde{K}^m_{\alpha\beta}(\epsilon) =
\tilde{e}_{\alpha\beta}(u^m(\epsilon)),~~
\tilde{K}^m_{\alpha 3}(\epsilon)
=\frac{1}{\epsilon}\tilde{e}_{\alpha 3}(u^m(\epsilon)),~~
\tilde{K}^m_{33}(\epsilon)
=\frac{1}{\epsilon^2}\tilde{e}_{33}(u^m(\epsilon)).\label{eq:e25}
\end{align}

Then there exists a constant $C_{10}>0$  such that
\begin{equation}
\|u^m(\epsilon)\|_{1,\Omega}\leq C_{10},~~~
|\tilde{K}^m_{ij}(\epsilon)|_{0,\Omega}\leq C_{10},~~~
|\tilde{\varphi}^m_i(\epsilon)|_{0,\Omega}\leq C_{10}\label{eq:e17}
\end{equation}
for all $0<\epsilon\leq\epsilon_0$.

Letting $v = u^m(\epsilon)$ in (\ref{eq:b13}), we have
\begin{align}
&\int_\Omega A^{ijkl}(\epsilon)
e_{k\|l}(\epsilon)(u^m(\epsilon))
e_{i\|j}(\epsilon)(u^m(\epsilon))\sqrt{g(\epsilon)}{\rm d}x\nonumber\\[.2pc]
&\quad\ +\int_\Omega P^{3kl}(\epsilon)\partial_3{\varphi}^m
(\epsilon)e_{k\|l}(\epsilon)(u^m(\epsilon))\sqrt{g(\epsilon)}{\rm d}x\nonumber\\[.2pc]
&\quad\ +\epsilon\int_\Omega P^{\alpha kl}(\epsilon)
\partial_\alpha{\varphi}^m(\epsilon)
e_{k\|l}(\epsilon)(u^m(\epsilon))\sqrt{g(\epsilon)}{\rm d}x\nonumber\\[.2pc]
&=\xi^m(\epsilon)\int_\Omega [\epsilon^2u^m_\alpha(\epsilon)
u^m_\alpha(\epsilon)+u^m_3(\epsilon)u^m_3(\epsilon)]\sqrt{g(\epsilon)}{\rm d}x.
\end{align}
Letting $\psi=\varphi^m(\epsilon)$ in (\ref{eq:b14}) and using it in
the above equation, we get
\begin{align}
&\int_\Omega A^{ijkl}(\epsilon) e_{k\|l}(\epsilon, u^m(\epsilon))
e_{i\|j}(\epsilon, u^m(\epsilon))\sqrt{g(\epsilon)}{\rm d}x\nonumber\\[.2pc]
&\quad\ + \int_\Omega {\cal E}^{ij}(\epsilon){\tilde{\varphi}}^m_{i}(\epsilon)
{\tilde{\varphi}}^m_{j}(\epsilon)\sqrt{g(\epsilon)}{\rm d}x\nonumber\\
&=\xi^m(\epsilon)\int_\Omega[\epsilon^2u^m_\alpha(\epsilon) \cdot
u^m_\alpha(\epsilon)+u^m_3(\epsilon)u^m_3(\epsilon)]\sqrt{g(\epsilon)}{\rm d}x.
\label{eq:e20}
\end{align}

Using the coerciveness properties (\ref{eq:c11}) and
(\ref{eq:b3}), the inequality
$(a-b)^2\geq a^2/2-b^2$ and the generalized Korn's inequality
(\ref{eq:c12}), we have for $\epsilon\leq \min\{\epsilon_0, 1\}$,
\begin{align}
&\int_\Omega A^{ijkl}(\epsilon) e_{k\|l}(\epsilon, u^m(\epsilon))
e_{i\|j}(\epsilon, u^m(\epsilon))\sqrt{g(\epsilon)}{\rm d}x\nonumber\\[.2pc]
&\quad\ + \int_\Omega {\cal E}^{ij}(\epsilon){\tilde{\varphi}}^m_{i}(\epsilon)
{\tilde{\varphi}}^m_{j}(\epsilon)\sqrt{g(\epsilon)}{\rm d}x\nonumber\\[.2pc]
&\geq  C_{11}\sum_{i,j}\|e_{i\|j}(\epsilon,
u^m(\epsilon))\|^2_{0,\Omega}
+C_{11}\sum_i\|\tilde{\varphi}^m_i(\epsilon)\|^2_{0,\Omega}\nonumber
\end{align}
\begin{align}
\hskip -1pc &= C_{11}\sum_{\alpha,\beta}
\|\tilde{e}_{\alpha\beta}(u^m(\epsilon))+\epsilon^2
{e}^\sharp_{\alpha\beta}(\epsilon,
u^m(\epsilon))\|^2_{0,\Omega}\nonumber\\[.2pc]
\hskip -1pc &\quad\ +2C_{11}\sum_{\alpha}
\left\|\frac{1}{\epsilon}\tilde{e}_{\alpha 3}(u^m(\epsilon))+\epsilon
{e}^\sharp_{\alpha 3}(\epsilon,
u^m(\epsilon))\right\|^2_{0,\Omega}\nonumber\\[.2pc]
\hskip -1pc &\quad\ +C_{11} \left\|\frac{1}{\epsilon^2}\tilde{e}_{33}(u^m(\epsilon))\right\|^2_{0,\Omega}
+C_{11}\sum_i \|\tilde{\varphi}^m_i(\epsilon)\|^2_{0,\Omega}\nonumber\\[.2pc]
\hskip -1pc &\geq C_{11}\left\lbrace \frac{1}{2}\sum_{i,j}|\tilde{K}_{ij}^m(\epsilon)|^2_{0,\Omega}
-C_1^2(2\epsilon^2+\epsilon^4)\|u^m(\epsilon)\|^2_{1,\Omega}\right\rbrace\nonumber\\[.2pc]
\hskip -1pc &\quad\ +C_{11}\sum_i\|\tilde{\varphi}^m_i(\epsilon)\|^2_{0,\Omega}\hskip 1pc \nonumber\\[.2pc]
\hskip -1pc &\geq C_{11}\left\lbrace \frac{1}{2}\sum_{i,j}
\|\tilde{e}_{ij}(u^m(\epsilon))\|^2_{0,\Omega} -
3\epsilon^2C_1^2\|u^m(\epsilon)\|^2_{1,\Omega}\right\rbrace\nonumber\\[.2pc]
\hskip -1pc &\quad\ +C_{11}\sum_i\|\tilde{\varphi}^m_i(\epsilon)\|^2_{0,\Omega}\nonumber\\[.2pc]
\hskip -1pc &\geq C_{11}\left\lbrace \frac{1}{2}(C_8)^{-2} -
3\epsilon^2C_1^2\right\rbrace \|u^m(\epsilon)\|^2_{1,\Omega}
+ C_{11}\sum_i\|\tilde{\varphi}^m_i(\epsilon)\|^2_{0,\Omega}.\label{eq:e21}
\end{align}

Combining eqs~(\ref{eq:e20}) and (\ref{eq:e21}) with relations
(\ref{eq:b15}) and (\ref{eq:d11}), we get the relation (\ref{eq:e17}).\vspace{.5pc}

\noindent {\it Step} (ii).\ \ From Step (i) it follows that there exists a
subsequence $(\tilde{\varphi}^m_i(\epsilon))$ and
$(\tilde{\varphi}^m_i)\in L^2(\Omega)$ such that
\begin{equation}
(\epsilon\partial_1\varphi^m(\epsilon),
\epsilon\partial_2\varphi^m(\epsilon),
\partial_3\varphi^m(\epsilon))\rightharpoonup (\tilde{\varphi}^m_1,
\tilde{\varphi}^m_2, \tilde{\varphi}^m_3)\quad \mbox{ in } (L^2(\Omega))^3.
\label{eq:e22}
\end{equation}
Since $\Gamma_{eD}$ contains $\Gamma^-$, we have
\begin{equation}
{\varphi}^m(\epsilon)(x_1, x_2, x_3)=
\int_{-1}^{x_3} \partial_3{\varphi}^m(\epsilon)(x_1, x_2, s){\rm d}s
\label{eq:e23}
\end{equation}
and it follows that $\|{\varphi}^m(\epsilon)\|_{0,\Omega} \leq
\sqrt{2}\|\partial_3{\varphi}^{m} (\epsilon)\|_{0,\Omega}$. This implies that
${\varphi}^m(\epsilon)$ is bounded in $L^2(\Omega)$. Therefore there
exists a $\varphi^m$ in $L^2(\Omega)$ and a subsequence, still indexed
by $\epsilon,$ such that $\varphi^m(\epsilon)$ converges weakly to
$\varphi^m$. Hence it follows from (\ref{eq:e22}) that
\begin{equation}
(\epsilon\partial_1\varphi^m(\epsilon),
\epsilon\partial_2\varphi^m(\epsilon),
\partial_3\varphi^m(\epsilon))\rightharpoonup
(0, 0, \partial_3\varphi^m).\label{eq:e24}
\end{equation}

\noindent {\it Step} (iii).\ \ From Step (i) it follows that there exists
a subsequence, indexed by $\epsilon$ for notational convenience,
and functions $u^m\in V(\Omega)$ and $\tilde{K}^m_{ij}\in
(L^2(\Omega))^9$ such that
\begin{equation}
u^m(\epsilon)\rightharpoonup u^m \quad\mbox{ in } H^1(\Omega),  ~~~
\tilde{K}^m(\epsilon)\rightharpoonup \tilde{K}^m\quad \mbox{ in }
L^2(\Omega)\mbox{ as }\epsilon\rightarrow 0.\label{eq:e28}
\end{equation}

Then there exist functions $(\zeta^m_\alpha)\in H^1(\omega)$ and
$\zeta^m_3\in H^2(\omega)$ satisfying $\zeta^m_i=\partial_\nu
\zeta^m_3=0$ on $\gamma_0$ such that
\begin{equation}
u^m_\alpha=\zeta^m_\alpha-x_3\partial_\alpha\zeta^m_3 \quad\mbox{ and }\quad
u^m_3=\zeta^m_3\label{eq:e29}
\end{equation}
and
\begin{align}
\tilde{K}^m_{\alpha\beta} &= \tilde{e}_{\alpha\beta}(u^m), ~~~
\tilde{K}^m_{\alpha 3}=-\frac{1}{\mu}P^{3\alpha 3}\partial_3\varphi^m,\nonumber\\[.2pc]
\tilde{K}^m_{33} &= -\frac{1}{\lambda+2\mu}(P^{333}\partial_3\varphi^m+\lambda
\tilde{K}^m_{\beta\beta}).\label{eq:e32}
\end{align}
From definition (\ref{eq:e25}) and the boundedness of
$(\tilde{K}^m_{ij}(\epsilon))$, we deduce that
\begin{equation*}
\|e_{\alpha 3}(u^m(\epsilon))\|_{0,\Omega}\leq \epsilon C_{13} \quad \hbox{and}\quad
\|e_{33}(u^m(\epsilon))\|_{0,\Omega}\leq\epsilon^2C_{13},
\end{equation*}
where $e_{ij}(v)=\frac{1}{2}(\partial_iv_j+\partial_jv_i)$.
Since norm is a weakly lower semicontinuous function
\begin{equation}
 \|e_{i3}(u^m)\|_{0,\Omega}\leq  \liminf_{\epsilon\rightarrow 0
}\|e_{i3}  (u^m(\epsilon)\|_{0,\Omega}=0,\label{eq:e30}
\end{equation}
we obtain $e_{i3}(u^m)=0$. Then it is a standard argument that the
components $u^m_i$ of the limit $u^m$ are of the form (\ref{eq:e29}).

Since $u^m(\epsilon)\rightharpoonup u^m$ in $H^1(\Omega)$,
definition (\ref{eq:c4})
of the functions $\tilde{e}_{\alpha\beta}(v)$ shows that the function
$\tilde{K}^m_{\alpha\beta}(\epsilon)=\tilde{e}_{\alpha\beta}(u^m(\epsilon))$
converges weakly in $L^2(\Omega)$ to the function
$\tilde{e}_{\alpha\beta}(u^m)$.

We next note the following result. Let $w\in L^2(\Omega)$ be given;
then
\begin{equation}
\int_\Omega w\partial_3v{\rm d}x=0 \quad\mbox{ for all }v\in H^1(\Omega) \mbox{
with } v=0\mbox{ on } \Gamma_0, \mbox{ then } w=0.
\end{equation}

Multiplying (\ref{eq:b13}) by $\epsilon^2$, taking $(v_\alpha)=0$
and letting  $\epsilon\rightarrow 0$, we get
\begin{equation}
\int_\Omega (\lambda{\tilde{K}^m}_{\sigma\sigma}
+(\lambda+2\mu)\tilde{K}_{33} + P^{333}\partial_3\varphi^m)
\partial_3v_3{\rm d}x =0 \label{eq:e33}
\end{equation}
which implies $(\lambda{\tilde{K}^m}_{\sigma\sigma}
+(\lambda+2\mu)\tilde{K}_{33}+P^{333}\partial_3\varphi^m)=0$
and hence the third relation in (\ref{eq:e32}) follows.

Again, multiplying (\ref{eq:b13}) by $\epsilon$, taking
$v_3=0$ and letting $\epsilon\rightarrow 0 $, we get
\begin{equation}
\int_\Omega (\mu{\tilde{K}^m}_{\alpha 3}+P^{3\alpha3}
\partial_3\varphi^m)\partial_3v_\alpha {\rm d}x = 0\label{eq:e34}
\end{equation}
which implies $(\mu{\tilde{K}^m}_{\alpha 3}+P^{3\alpha3}
\partial_3\varphi^m)=0$ and hence the second relation in (\ref{eq:e32})
follows.\vspace{.5pc}

\noindent{\it Step} (iv).\ \ The function $\varphi^m$ is of the form
(\ref{eq:e10}).

Letting $\epsilon\rightarrow 0$ in eq.~(\ref{eq:b14}), we get
\begin{equation}
\int_\Omega (P^{3\alpha\beta}\tilde{K}^m_{\alpha\beta} -
{\cal E}^{33}\partial_3\varphi^m)\partial_3\psi {\rm d} x = 0\quad
\forall \psi\in \Psi(\Omega).\label{eq:e38}
\end{equation}

Since $D(\Omega)$ is dense in $\Psi_{l0}$ (and hence in $\Psi(\Omega)$)
for the norm $\|.\|_{\Psi_l}$, eq.~(\ref{eq:e38}) is equivalent to
\begin{equation}
\partial_3(P^{3\alpha\beta}\tilde{K}^m_{\alpha\beta} -
{\cal E}^{33}\partial_3\varphi^m) = 0\quad \mbox{ in } D'(\Omega)\label{eq:e39}
\end{equation}
which implies that $(P^{3\alpha\beta}\tilde{K}^m_{\alpha\beta}-
{\cal E}^{33}\partial_3\varphi^m)=d^1$, with $d^1\in D(\omega)$. Then
\begin{equation}
\partial_3\varphi^m=\frac{p^{3\alpha\beta}}{p^{33}}
[\tilde{e}_{\alpha\beta}(\zeta^m)
-x_3\partial_{\alpha\beta}\zeta^m_3]-\frac{1}{p^{33}}d^1\label{eq:e40}
\end{equation}
which gives
\begin{equation}
\varphi^m=\frac{p^{3\alpha\beta}}{p^{33}}[x_3\tilde{e}_{\alpha\beta}(\zeta^m)
-x_3^2\partial_{\alpha\beta}\zeta^m_3]-\frac{x_3}{p^{33}}d^1+d^0.
\label{eq:e41}
\end{equation}

Since $\varphi^m$ satisfies the boundary conditions
$\varphi^m_{|\Gamma^+} =\varphi^m_{|\Gamma^-}=0$, we have
\begin{equation}
d^0 = \frac{p^{3\alpha\beta}}{2p^{33}}
\partial_{\alpha\beta}\zeta^m_3, ~~~~
d^1=p^{3\alpha\beta}\tilde{e}_{\alpha\beta}(\zeta^m).\label{eq:e43}
\end{equation}
Thus the conclusion follows.\vspace{.5pc}

\noindent {\it Step} (v).\ \ The function $(\zeta^m_i)$ satisfies
(\ref{eq:e11}) and (\ref{eq:e12}).

Taking $v\in V_{KL}$ and letting $\epsilon\rightarrow 0$ in
(\ref{eq:b13}) we get
\begin{equation}
\int_\Omega A^{\alpha\beta kl}\tilde{K}^m_{kl}\tilde{K}_{\alpha\beta}(v){\rm d}x
+\int_\Omega P^{3\alpha\beta}\partial_3\varphi^m \tilde{K}_{\alpha\beta}(v){\rm d}x
=\xi^m\int_\Omega u^m_3 \cdot v_3 {\rm d}x.\label{eq:e44}
\end{equation}
Replacing $u^m$ and $\tilde{K}^m_{ij}$ by the expressions obtained in
(\ref{eq:e29}) and (\ref{eq:e32}), and taking $v$ of the form
\begin{equation*}
v_\alpha=\eta_\alpha-x_3\partial_\alpha\eta_3\quad \mbox{ and }\quad
v_3=\eta_3
\end{equation*}
with $(\eta_i)\in V_H(\omega)\times V_3(\omega)$,
it is verified that (\ref{eq:e44}) coincides with
eqs~(\ref{eq:e11}) and (\ref{eq:e12}).\vspace{.5pc}

\noindent {\it Step} (vi).\ \ The convergences $u^m(\epsilon)\rightharpoonup u^m$ in
$H^1(\Omega)$ and $\varphi^m(\epsilon)\rightharpoonup \varphi^m$
in $L^2(\Omega)$ are strong.

To show that the family $(u^m(\epsilon))$ converges strongly to $u^m$
in $H^1(\Omega),$ by Lemma~4.2, it is enough to show that
\begin{equation}
\tilde{e}_{ij}(u^m(\epsilon))\rightarrow \tilde{e}_{ij}(u^m)\quad\mbox{ in }L^2(\Omega).
\label{eq:e46}
\end{equation}
Since $\tilde{e}_{i3}(u^m)=0$ and
\begin{align}
&\sum_{i,j}\|\tilde{e}_{ij}(u^m(\epsilon))-\tilde{e}_{ij}(u^m)\|^2_{0,\Omega}\nonumber\\[.2pc]
&=\sum_{\alpha,\beta}\|\tilde{K}^m_{\alpha\beta}(\epsilon)-
\tilde{K}^m_{\alpha\beta}\| _{0,\Omega}^2 +2\epsilon^2\sum_{\alpha}\|\tilde{K}^m_{\alpha
3}(\epsilon)\|^2_{0,\Omega}+\epsilon^4\|\tilde{K}^m_{33}(\epsilon)\|^2_{0,\Omega},\label{eq:e47}
\end{align}
convergence (\ref{eq:e46}) is equivalent to showing that
\begin{equation}
\tilde{K}^m(\epsilon)\rightarrow \tilde{K}^m \quad\mbox{ in } L^2(\Omega).
\label{eq:e48}
\end{equation}

We define a norm on $(L^2(\Omega))^9\times (L^2(\Omega))^3$ by letting
for any matrix $M\in (L^2(\Omega))^9$ and any vector $\chi\in
(L^2(\Omega))^3$,
\begin{equation}
\|(M, \chi)\| = \left\lbrace \int_\Omega A^{ijkl}M:M\sqrt{g(\epsilon)}{\rm d}x
+\int_\Omega{\cal E}^{ij}\chi_i\chi_j \sqrt{g(\epsilon)}{\rm d}x\right\rbrace^{1/2}.\label{eq:e49}
\end{equation}

Let $X^m(\epsilon)$ be the norm of $(\tilde{K}^m(\epsilon),
\epsilon\partial_1\varphi^m(\epsilon),
\epsilon\partial_2\varphi^m(\epsilon), \partial_3\varphi^m(\epsilon)) $
in $(L^2(\Omega))^{12}$. Using the weak convergence equation (eqs~(\ref{eq:e24}) and
(\ref{eq:e28})) and the relation (\ref{eq:e32}), it can be shown that
\begin{equation}
\lim_{\epsilon\rightarrow 0}X^m(\epsilon)=X^m= \left( \int_\Omega
A^{ijkl}\tilde{K}^m :\tilde{K}^m{\rm d}x +\int_\Omega {\cal E}^{33}(\partial_3\varphi^m)^2
{\rm d}x \right)^{1/2}\label{eq:e50}
\end{equation}
which is the norm of $(\tilde{K}^m, 0, 0, \partial_3\varphi^m)$. Since
we have already proved that
$(\tilde{K}^m(\epsilon)$, $\epsilon\partial_1\varphi^m(\epsilon)$,
$\epsilon\partial_2\varphi^m(\epsilon)$, $\partial_3\varphi^m(\epsilon))$
converges weakly to $(\tilde{K}, 0, 0, \partial_3\varphi^m)$ in
$(L^2(\Omega))^{12}$, we have the following strong convergences:
\begin{align}
&\tilde{K}^m(\epsilon)\rightarrow \tilde{K}^m \mbox{ strongly in }
(L^2(\Omega))^9,\label{eq:e51}\\[.2pc]
&(\epsilon\partial_1\varphi^m(\epsilon),
\epsilon\partial_2\varphi^m(\epsilon), \partial_3\varphi^m(\epsilon))
\rightarrow (0,0,\partial_3\varphi^m)\mbox{ strongly in }
(L^2(\Omega))^3.\label{eq:e52}
\end{align}

Hence $u^m(\epsilon)$ converges strongly to $u^m$ in $H^1(\Omega)$ and
since $\varphi^m(\epsilon)-\varphi^m$ is in $\Psi_{l0}$, the equivalence
of norms $\|\psi\|_{\Psi_l}$ and $\psi\rightarrow
|\partial_3\psi|_\Omega$ in $\Psi_{l0}$ proves that
$\varphi^m(\epsilon)$ converges strongly to $\varphi^m$ in
$L^2(\Omega)$.\hfill \ab
\end{proof}

Equation (\ref{eq:e12}) can be written as
\begin{align}
&\int_\omega \left[\frac{2\lambda\mu}{\lambda+2\mu}e_{\rho\rho}(\zeta)
\delta_{\alpha\beta}+2\mu e_{\alpha\beta}(\zeta) \right]
\partial_\beta\eta_\alpha {\rm d}\omega\nonumber\\[.2pc]
&= \int_\omega \left[ \frac{2\lambda\mu}{\lambda+2\mu}(\partial_\sigma\theta
\partial_\sigma\zeta_3)\delta_{\alpha\beta}
+\mu(\partial_\alpha\theta\partial_\beta\zeta_3+\partial_\beta\theta
\partial_\alpha\zeta_3) \right]\partial_\beta\eta_\alpha {\rm d}\omega.\label{eq:e53}
\end{align}
Clearly, the  bilinear form
\begin{align}
\tilde{b}(\zeta_\alpha, \eta_\alpha) &=
\int_\omega \left[\frac{2\lambda\mu}{\lambda+2\mu}e_{\rho\rho}(\zeta)
\delta_{\alpha\beta}+2\mu e_{\alpha\beta}(\zeta) \right]
\partial_\beta\eta_\alpha {\rm d}\omega\nonumber\\[.2pc]
&= \int_\omega \left[\frac{2\lambda\mu}{\lambda+2\mu}e_\rho\rho(\zeta)
e_{\sigma\sigma}(\eta)+2\mu e_{\alpha\beta}(\zeta)e_{\alpha\beta}(\eta)\right]{\rm d}\omega
\label{eq:e54}
\end{align}
is $V_H(\omega)$ elliptic. Also for a  given $\zeta_3\in V_3(\omega)$,
the functional

\begin{equation}
\langle \zeta_3, \eta_\alpha \rangle = \int_\omega
\left[ \frac{2\lambda\mu}{\lambda+2\mu}(\partial_\sigma\theta
\partial_\sigma\zeta_3)\delta_{\alpha\beta}
+\mu(\partial_\alpha\theta\partial_\beta\zeta_3+\partial_\beta\theta
\partial_\alpha\zeta_3) \right]\partial_\beta\eta_\alpha {\rm d}\omega
\end{equation}
is continous on $V_H(\omega)$. Thus, given $\zeta_3 \in V_3(\omega)$,
there exists a unique vector $(\zeta_\alpha)\in V_H(\omega)$ such that
\begin{equation}
\tilde{b}(\zeta_\alpha, \eta_\alpha) = \langle \zeta_3, \eta_\alpha\rangle.
\end{equation}
We denote by $T\zeta_3\in V_H(\omega)\times V_3(\omega)$ the vector
$(\zeta_\alpha, \zeta_3)$. In particular, $T\zeta^m_3=(\zeta^m_\alpha,
\zeta^m_3)$.

Substituting this in (\ref{eq:e11}), we get
\begin{equation}
b(\zeta^m_3, \eta_3)=\xi^m\int_\omega \zeta^m\eta_3 {\rm d}\omega
\quad \mbox{ for all } \eta_3\in V_3(\omega),\label{eq:e55}
\end{equation}
where
\begin{align}
b(\zeta_3, \eta_3) &= -\int_\omega m_{\alpha\beta}\partial_{\alpha\beta}\eta_3 {\rm d}\omega +
\int_\omega n^\theta_{\alpha\beta}(T\zeta_3)\partial_{\alpha\beta}\theta\eta_3
{\rm d}\omega\nonumber\\[.2pc]
&\quad +\frac{2}{3}\int_\omega \frac{p^{3\alpha\beta}p^{3\rho\tau}}{p^{33}}
\partial_{\rho\tau}\zeta_3\partial_{\alpha\beta}\eta_3 {\rm d}\omega.
\label{eq:e56}
\end{align}

\begin{lem}
The bilinear form $b(\cdots)$ defined by {\rm (\ref{eq:e56})} is
$V_H(\omega)$-elliptic and symmetric.
\end{lem}

\begin{proof}
It follows from Lemma~6.2 in \cite{SKNS1} that the bilinear form
$\tilde{b}(\cdots)$ defined by
\begin{equation}
\tilde{b}(\zeta_3, \eta_3)=
-\int_\omega m_{\alpha\beta}(\zeta_3)\partial_{\alpha\beta}\eta_3
{\rm d}\omega + \int_\omega
n^\theta_{\alpha\beta}(T\zeta_3)\partial_{\alpha\beta}\theta\eta_3
{\rm d}\omega\label{eq:e57}
\end{equation}
is $V_H(\omega)$-elliptic and symmetric. Hence it is clear that $b(\cdots)$ is also
$V_H(\omega)$-elliptic and symmetric.\hfill \ab
\end{proof}

\begin{lem}
Let $(\zeta_3^m, \xi^m), m\geq 1${\rm ,} be the eigensolutions of problem
{\rm (\ref{eq:e56})} found as limits of the subsequence $(u^m(\epsilon){\rm ,}
\xi^m(\epsilon)), m\geq 1$ of eigensolutions of the problem
{\rm (\ref{eq:b13})}. Then the sequence $(\xi^m)_{m=1}^{\infty}$ comprises all
the eigenvalues{\rm ,} counting multiplicities{\rm ,} of problem {\rm (\ref{eq:e56})} and
the associated sequence $(\zeta^m_3)_{m=1}^{\infty}$ of eigenfunctions
forms a complete orthonormal set in the space $V_3(\omega)$.
\end{lem}

\begin{proof}
The proof is similar to the proof of Lemma~5.4 in \cite{CK}.\hfill \ab
\end{proof}

\end{document}